\documentclass[a4paper,11pt]{article}
\usepackage{amsfonts}
\usepackage{amsmath}
\usepackage{amssymb}

\usepackage[latin1]{inputenc}

\newcommand{\R}{\mathbb{R}}

\newcommand{\N}{\mathbb{N}}

\newcommand{\proof}[1]{\par\smallskip\noindent{\bf Proof#1.}}

\newcommand{\ve}{\varepsilon}

\def \e {\varepsilon}
\def \ue {u^{\varepsilon}}
\def \ve {v^{\varepsilon}}
\def \me {m^{\varepsilon}}
\def \he {h^{\varepsilon}}
\def \p{\partial}
\def \O{\Omega}

\newtheorem{lem}    {Lemma}[section]
\newtheorem{pro}    [lem]{Proposition}
\newtheorem{thm}    [lem]{Theorem}

\newtheorem{cor}    [lem]{Corollary}
\newtheorem{df}     [lem]{Definition}

\renewcommand\theequation{\arabic{section}.\arabic{equation}}
\title{The singular limit of a haptotaxis model with bistable growth}
\author{Elisabeth LOGAK $^{\ast}$  and  Chao WANG \!
\thanks {University of Cergy-Pontoise, Department of Mathematics, UMR CNRS 8088,
F-95000 Cergy-Pontoise}}
\date{}
\begin{document}
\maketitle
\begin{abstract}
\noindent
We consider a model for haptotaxis with bistable growth and study its singular limit.
This yiels an interface motion where the normal velocity of the interface depends
on the mean curvature and on some nonlocal haptotaxis term. We prove the result for general
initial data after establishing a result about generation of interface in a small time.
\end{abstract}
\noindent

\def\theequation{\arabic{section}.\arabic{equation}}

\section{Introduction}
\setcounter{equation}{0}
In this article we consider a model for haptotaxis with growth.
Haptotaxis is the directed motion of cells by migration up a gradient
of cellular adhesion sites located in the extracellular matrix (ECM).
This process appears in tumor invasion and is involved in the first stage
 of proliferation. It also plays an important role in wound healing. \\

\noindent 
The basic mechanism involves 3 main cellular components: the tumor cells,
the Extracellular Matrix (ECM), and some Matrix Degrading Enzymes (MDE).
Tumor cells migrate in response to gradients of some ECM proteins. Those ECM proteins
are degraded by MDE, those enzymes being produced by tumor cells themselves.
 Moreover, both tumor cells and MDE diffuse in the cellular medium but
 ECM proteins do not diffuse.
 
\noindent
This mechanism is reminiscent of chemotaxis, which  is accounting
for the directed migration of biological individuals
 (e.g. bacteria) towards higher gradients of some chemical substance.
Chemotaxis often works as an aggregating mechanism, which
is reflected in the blow-up of solutions of the Keller-Segel model, a phenomenon that
 has been widely studied in the recent years.
However there is a major difference between chemotaxis and haptotaxis:
 since ECM proteins do not diffuse, instead of the elliptic or parabolic coupling appearing in chemotaxis,
the haptotaxis model involves an ODE coupling between the
concentration of ECM proteins and the MDE concentration. This is also the case in angiogenesis model
(cf \cite{cpz}) but models for haptotaxis  involve (at least) 3 equations, whereas angiogenesis is a coupled system of 2 equations.\\

\noindent
We now give a brief review of the mathematical literature related to haptotaxis modelling.
The relevant variables are the tumor cells concentration, the Extracellular
Matrix concentration (ECM) the Matrix Degrading Enzymes concentration (MDE)
 as well as the oxygen concentration. A hybrid model using PDEs and
cellular automata has been proposed by Anderson
 \cite{anderson05},   involving
4 components: tumor cells, ECM, MDE + Oxygen.
 Global existence for Anderson's model has been established
in \cite{ww} (Walker, Webb (07)).
Our model  is a simpler version from this model involving 3 components
where we introduce a
bistable nonlinearity to model the role of
 changes in oxygen concentration.
A similar model of haptotaxis with a logistic nonlinearity is studied in \cite{pmc09} 
(and the references therein)
and global existence is proved.
Finally Chaplain, Lolas (see \cite{cl06} and the references therein) proposed a combined
 chemotaxis-haptotaxis model with logistic source.
 Recent results by Y. Tao, M. Wang \cite{tw08} and  Y. Tao \cite{tw09} show
 global existence for this model in dimension $N\leq 2$.
Complex patterns in haptotaxis models
are obtained numerically in \cite{ww}, and also in \cite{yekat} and
\cite{chertock}.

\noindent
Our starting point is the haptotaxis model proposed in \cite{ww}. In
this paper, the authors prove global well-posedness for a large class of initial
data, a result which strongly emphasizes the difference with Keller-Segel chemotaxis
model.  Here we consider a different version of this model, where we do not explicitely
consider the oxygen concentration as a variable. Instead we replace it by a bistable
nonlinearity in the equation for the cell concentration.\\
Next we show that in the limit $\varepsilon \rightarrow 0$, the solutions converge to the
solutions of a free boundary problem where the interface motion is driven by mean
curvature plus an haptotaxis term.

\noindent
More precisely, we study the initial value Problem $(P^{\varepsilon})$
$$
  (P^{\varepsilon})
 \left\{
\begin{array}{ll}
u_{t} = \Delta u - \nabla \cdot (u \nabla \chi (v)) +
 \dfrac{1}{\varepsilon^{2}} f(u) & \quad \quad \mbox{ in } \ \O \times (0,T]  \\
v_{t}= -\lambda mv & \quad \quad \mbox{ in }  
\ \O \times (0,T] \\
m_{t} = \alpha \Delta m +u-m & \quad \quad \mbox{in}  \ \O \times (0,T] \\
u(x,0)=u_{0}(x) & \quad \quad x \in \Omega \\
v(x,0)=v_{0}(x) & \quad \quad x \in \Omega \\
m(x,0)=m_{0}(x) & \quad \quad x \in \Omega \\
\dfrac{\p u}{\p \nu} = \dfrac{\p m}{\p \nu} =0  & \quad \quad \mbox{on} \ \p \O \times (0,T],
\end{array}
\right.
$$
where $\O$ is a smooth bounded domain in $\mathbb{R}^{N}$ $(N \geq 2)$,
$\O_T=\Omega \times [0,T]$ with $T >0$,
$\nu$ is the exterior normal vector on $\p \O$  and  $\lambda >0$, $\alpha >0$
are strictly positive constants.

\noindent
 The haptotaxis sensitivity function $\chi$ is smooth and satisfies
 $$\forall v>0,\,\,\,\chi(v) >0, \,\,\,\chi'(v)>0 .$$
The growth term $f$ is bistable and is given by
$$\forall u \in \R,\,\,\,f(u) = u(1-u)(u-\frac{1}{2}) $$
so that $\int_{0}^{1}f(u)du=0$.

\noindent
We make the following assumptions about the initial data.
\begin{enumerate}
\item $u_{0}$, $v_{0}$ and $m_{0}$ are nonnegative $C^{2}$ functions in
 $\overline \O$ and we fix a constant $C_{0} >1$ such that
\begin{equation}\label{diyi}
 ||u_{0}||_{C^{2}(\overline \O)} +|| v_{0}||_{C^{2}(\overline \O)}
+||m_{0}||_{C^{2}(\overline \O)} \leq C_{0}.
\end{equation}
\item
$v_{0}$ satisfies the homogeneous Neumann boundary condition
\begin{equation}\label{neumann0}
 \dfrac{\p v_{0}}{\p \nu} =0 \quad \mbox{on} \quad \p \O.
\end{equation}
\item
The open set $\O_{0}$ defined by
$$\O_{0} := \{ x \in \O, u_{0}(x) > 1/2 \}$$
 is connected and $\O_{0} \subset \subset \O$.
\item
$\Gamma_{0} := \p \O_{0}$ is a smooth hypersurface without boundary.
\end{enumerate}
With these assumptions $\O_{0}$ is a domain enclosed by the initial interface $\Gamma_{0}$ and
$$u_{0} > 1/2 \mbox{ in }  \O_{0},\,\,\, 0 \leq u_{0} < 1/2
\mbox{ in }  \O \setminus \overline \O_{0}.$$
The existence of a unique nonnegative solution $(u^{\varepsilon}, v^{\varepsilon}, m^{\varepsilon})$
 to Problem $(P^{\varepsilon})$ is established in Section 2.  Note that it follows
from (\ref{neumann0}) that
\begin{equation}\label{neumann}
 \frac{\p v^{\varepsilon}}{\p \nu} =0 \mbox{ on }  \p \O \times [0,T].
\end{equation}
  We are interested
 in the asymptotic behavior of $(u^{\varepsilon}, v^{\varepsilon}, m^{\varepsilon})$
as $\varepsilon \rightarrow 0$.
The asymptotic limit of Problem $(P^{\varepsilon})$
as $\varepsilon \rightarrow 0$ is given by the following free boundary
Problem $(P^{0})$
\begin{align*}
 (P^{0})
 \left\{
\begin{array}{ll}
u^{0}(x,t)=  \chi_{\O_{t}}(x) = \displaystyle{
\left\{
\begin{array}{ll}
1  \mbox{   in } \ \O_{t}, t \in [0,T]   \\
0  \mbox{   in } \ \O \setminus \overline \O_{t}, t \in [0,T]
\end{array}
\right.}  \\
v_{t}^{0} =-\lambda m^{0}v^{0}  & \mbox{in}  \ \O \times (0,T] \\
m_{t}^{0}= \alpha \Delta m^{0} + u^{0} - m^{0} &  \mbox{in}  \ \O \times (0,T] \\
V_{n} = -(N-1) \kappa + \dfrac{\p \chi (v^{0})}{\p n} &   \mbox{on} \ \Gamma_{t}= \p \O_{t}, t \in (0, T]  \\
\Gamma_{t} |_{t=0} = \Gamma_{0} \\
v^{0}(x,0)=v_{0}(x) &   x \in \Omega \\
m^{0}(x,0)=m_{0}(x)  & x \in \Omega \\
\dfrac{\p m^{0}}{\p \nu} =0 &    \mbox{on} \ \p \O \times (0,T],
\end{array}
\right.
\end{align*}
\noindent
where $\O_{t} \subset \subset \O$ is a moving domain,
$\Gamma_{t}= \p \O_{t}$ is the limit interface,
$n$ is the exterior normal vector on $\Gamma_{t}$ $V_{n}$ is the normal velocity of $\Gamma_{t}$
in the exterior direction and $\kappa$ is the mean curvature at each point
of $\Gamma_{t}$.
 We first establish the well-posedness of Problem $(P^{0})$ locally in time in Section 3. 
Our main result is to prove rigorously the convergence of
$(u^{\varepsilon}, v^{\varepsilon}, m^{\varepsilon})$
to
$(u^0, v^0, m^0)$
for initial data satisfying the above assumptions.
In a first step, we establish the following generation of interface property.
\begin{thm}
\label{gmi}
Assume that $(u_{0},v_0,m_0)$ satisfy the hypotheses $1$-$2$-$3$-$4$.
 Let $0 <\eta < 1/4$ be an arbitrary constant and define $\mu = f^{'}(1/2) = 1/4$. Then there exist $\varepsilon_{0} >0$ and $M_{0}>0$ such that, for all $\varepsilon \in (0, \varepsilon_{0}]$ and all $t \in [t^{\ast},T]$ where $t^{\ast}=\mu^{-1} \varepsilon^{2} |\ln \varepsilon|$,\\
(a) for all $x \in \O$, we have
$$0 \leq u^{\varepsilon} (x, t) \leq 1+\eta ;$$
(b) for all $x \in \O$ such that $|u_{0}(x) - \frac{1}{2}| \geq M_{0} \varepsilon$, we have
$$\mbox{ if } u_{0}(x) \geq \frac{1}{2} + M_{0} \varepsilon , \mbox{ then }
 u^{\varepsilon}(x, t) \geq 1-\eta ,$$
$$\mbox{ if }  u_{0}(x) \leq \frac{1}{2} - M_{0} \varepsilon ,  \mbox{ then }
0\leq u^{\varepsilon}(x, t) \leq \eta .$$
\end{thm}
\noindent
The main result reads as follows.
\begin{thm}
\label{thm1}
Assume that $(u_{0},v_0,m_0)$ satisfy the hypotheses $1$-$2$-$3$-$4$. 
Let $(u^{\varepsilon}, v^{\varepsilon}, m^{\varepsilon})$ be the solution 
of Problem $(P^{\varepsilon})$ and let $(v^{0}, m^{0}, \Gamma)$ 
with $\Gamma = (\Gamma_{t} \times \{ t \}) _{t \in [0,T]}$ 
be the smooth solution of the free boundary Problem $(P^{0})$ on $[0,T]$. 
Then, as $\varepsilon \rightarrow 0$, the solution 
$(u^{\varepsilon}, v^{\varepsilon}, m^{\varepsilon})$ converges to
 $(u^{0}, v^{0}, m^{0})$ almost everywhere
 in $\bigcup_{0 < t \leq T} ((\O \setminus \Gamma_{t}) \times {t})$. 
More precisely,
$$\lim_{\varepsilon \rightarrow 0} u^{\varepsilon}(x,t) = u^{0}(x,t) 
\mbox{  a.e.  } \mbox{ in } 
\bigcup_{0 < t \leq T} ((\O \setminus \Gamma_{t}) \times {t}), $$
and for all $\alpha \in (0,1)$,
 $$\lim_{\varepsilon \rightarrow 0} 
||v^{\varepsilon} - v^{0}||_{C^{1+\alpha, (1+\alpha)/2}(\overline \O_{T})} =0,$$
$$\lim_{\varepsilon \rightarrow 0} 
||m^{\varepsilon}-m^{0}||_{C^{1+\alpha, (1+\alpha)/2}(\overline \O_{T})} =0.$$
\end{thm}
 We actually prove a stronger convergence result concerning $\ue$.
\begin{cor}
\label{coru}
Assume that $(u_{0},v_0,m_0)$ satisfy the hypotheses $1$-$2$-$3$-$4$.
Then for any $ t \in (0,  T]$,
\begin{eqnarray}
&& \lim_{\varepsilon \rightarrow 0} \ue(x,t)=  \chi_{\O_t}(x) =\displaystyle{
\left\{
\begin{array}{ll}
 1 & \mbox{ for } x \in \Omega_{t} \\
 0 & \mbox{ for } x \in \O \setminus {\overline{\Omega_{t}}}
\end{array}
\right.}
\end{eqnarray}
\end{cor}
Moreover like in \cite{a}, we also obtain the following estimate
of the distance between the interface $\Gamma_{t}$ solution of Problem $(P^0)$
and the set
$$\Gamma_{t}^{\varepsilon}:=\{ x \in \O, u^{\varepsilon}(x,t)=1/2 \}.$$
\begin{thm}
\label{thm2}
There exists $C>0$ such that
$$\Gamma_{t}^{\varepsilon} \subset \aleph_{C\varepsilon}(\Gamma_{t}) \mbox{ for } 0 \leq t \leq T, $$
where $\aleph_{r}(\Gamma_{t}):= \{ x \in \O, dist(x, \Gamma_{t})<r  \}$ is the
 tubular neighborhood of $\Gamma_{t}$ of radius $r>0$.
\end{thm}
\begin{cor}
\label{haus}
$\Gamma_{t}^{\varepsilon} \rightarrow \Gamma_{t}$ as $\varepsilon \rightarrow 0$,
 uniformly in $t \in [0,T]$ in the sense of the Hausdorff distance.
\end{cor}
The organization of the paper is as follows.
In section 2 we prove some a priori estimates,
establish a comparison principle for Problem $(P^{\varepsilon})$ and prove the
 existence of a unique global solution.
 In section 3 we prove the well-posedness of the free boundary problem $(P^{0})$
 and obtain the existence of a smooth unique solution up to some time $T>0$.
In section 4 we establish the property of generation of interface. Finally in section 5
 we prove the convergence of the solution
 of Problem $(P^{\varepsilon})$ to the solution of Problem $(P^{0})$.

\section{A priori estimates and comparison principle}
\setcounter{equation}{0}
\subsection{A priori estimates}
For a given $T>0$ and a given nonnegative function $u_{0} \in C^{2}(\overline \O)$,
we define
$$X_{T} = \{   u \in C^{0}( \overline \O_{T}),\quad 0 \leq u \leq C_{0} \mbox{ in } \O_{T} \mbox{ and } u(x,0)=u_{0}(x)  \},$$
where $C_{0} >1$ is the constant defined in \eqref{diyi}.
\noindent
It is convenient to rewrite Problem $(P^{\varepsilon})$ as an evolution equation
for $u$ with a nonlocal coefficient $H(u)=v$, namely
\begin{equation}
\left\{
\begin{array}{ll}
u_{t} = \Delta u - \nabla \cdot (u \nabla \chi(H(u)) ) +\dfrac{1}{\varepsilon^{2}} f(u) & \quad \quad \mbox{in} \ \O \times (0,T]  \\
u(x,0)=u_{0}(x), & \quad \quad x \in \Omega \\
\dfrac{\p u}{\p \nu} =0  & \quad \quad \mbox{on} \ \p \O \times (0,T],
\end{array}
\right.
\end{equation}
where for a given function $u=u(x,t) \in X_{T}$, we define $H(u) = v$
as the first component of the
unique solution $(v,m)$ of the auxiliary problem
\begin{equation}
\label{yeshi}
\left\{
\begin{array}{ll}
v_{t}= -\lambda m v & \quad \quad \mbox{in  } \Omega \times (0,T]\\
m_{t}= \alpha \Delta m + u-m & \quad \quad \mbox{in  } \Omega \times (0,T]\\
v(x,0)=v_{0}(x), & \quad \quad x \in \Omega \\
m(x,0)=m_{0}(x), & \quad \quad x \in \Omega \\
\dfrac{\p m}{\p \nu} =0  & \quad \quad \mbox{on  } \p \O \times (0,T].
\end{array}
\right.
\end{equation}
The functions $v_{0}$ and $m_{0}$ are given and satisfy $1$-$2$.
\noindent
We give below some a priori estimates on the solution to Problem $(P^\e)$ and state the
related properties of $H$.
\begin{lem}
\label{xinjia}
For $u \in X_{T}$, let $(v, m)$ be the solution of Problem \eqref{yeshi} and let
$H : X_{T} \rightarrow C^{2}( \overline \O_{T})$ be the operator defined by
 $H(u)=v$.
 Then there exists $C >0$ only depending on $T$ and $\O$ such that
\begin{enumerate}
 \item[](a) for all $(u_{1}, u_{2}) \in X_{T}^{2}$ with
 $0 \leq u_{1} \leq u_{2}$ in $\O_{T}$, the solution $(v_i,m_i)$ of
 Problem \eqref{yeshi} for $i=1,2$ satisfies
$$0 \leq m_{1} \leq m_{2}   \mbox{ and }\leq v_{2} \leq v_{1} \mbox{ in } \O_{T}$$
 so that the operator $H$ is nonincreasing on $X_T$.
\item[](b) for all $u \in X_{T}$,
$$ ||m||_{C^{1+\alpha,(1+\alpha)/2}(\overline \O_{T})}  \leq C C_{0}
 \mbox{ and }
 \sup_{(x,t) \in \overline \O_{T}} \big{|}\int^{t}_{0} \Delta m(x,s)ds \big{|}
\leq C C_{0}. $$
\item[](c) for all $u \in X_{T}$, the function $v = H(u)$ satisfies
$$||v||_{C^{0}(\overline \O_{T})} \leq C_{0} \quad
\mbox{ and } \quad ||\nabla v||_{C^{0}(\overline \O_{T})}
+ ||\Delta v||_{C^{0}(\overline \O_{T})} \leq CC_{0}^{3}. $$
\end{enumerate}
\end{lem}
\proof{}
To prove property $(a)$, let $ (u_{1}, u_{2}) \in X_{T}^{2}$ with
 $0 \leq u_{1} \leq u_{2}$ in $\O_{T}$. Since for $i = 1, 2$
$$(m_{i})_{t} - \alpha \Delta m_{i} + m_{i} = u_{i}\geq 0 \mbox{  in  } \O_{T},  $$
with
 $$m_{i}|_{t=0} = m_{0}\geq 0 \mbox{ and } \dfrac{\p m_i}{\p \nu} =0 \mbox{ on }
 \p \O \times (0,T],$$
we deduce from the standard maximum principle that
 $0 \leq m_{1} \leq m_{2}$ in $\O_T$.

\noindent
 Next solving the equation $v_{t}= -\lambda m v$ we get that
\begin{align}
\begin{split}
\label{guanyuv}
v_{i}(x,t) = v_{0}(x) e^{-\lambda \int_{0}^{t} m_{i}(x,s) ds}
\end{split}
\end{align}
for all $(x,t) \in \O_{T}$ and $i = 1, 2$,
so that $v_{1} \geq v_{2} \geq 0 $ in $\O_{T}$, which proves that $H$
is nonincreasing on $X_T$.

\noindent
In order to prove $(b)$ and $(c)$, note that $m$ satisfies
 the linear parabolic equation
\begin{equation}
\left\{
\begin{array}{ll}
m_{t} = \alpha \Delta m +u-m & \mbox{ in }  \ \O \times (0,T] \\
m(x,0)=m_{0}(x) & x \in \Omega \\
\dfrac{\p m}{\p \nu} =0  & \mbox{ on } \ \p \O \times (0,T]
\end{array}
\right.
\end{equation}
with $0\leq u \leq C_0$ in $\O_T$ and $ m_{0} \geq 0 \in \Omega $. Thus
 it follows from the maximum principle and from standard parabolic estimates that
there exists a constant
$C>0$ only depending on $T$ and $\O$ such that
\begin{equation}\label{estm}
 0\leq m \leq C_0 \mbox{ in } \O_{T}, \,\,\,\,\,
||m||_{C^{1+\alpha,(1+\alpha)/2}(\overline \O_{T})}  \leq C C_{0}.
\end{equation}
In view of \eqref{guanyuv}, $v \geq 0$ and $v_{t} \leq 0$ in $\O_{T}$
 so that
\begin{align}
\begin{split}
0 \leq v(x,t) \leq v_{0}(x) \leq C_{0} \mbox{ for all } (x, t) \in \O_{T}.
\end{split}
\end{align}
Since for all $(x,t) \in \O_T$
\begin{align}
\begin{split}
\label{erdba}
\nabla v(x,t)= 
\nabla v_{0}(x) e^{-\lambda \int_{0}^{t} m(x,s) ds} -
 \lambda v(x,t) \big(\int_{0}^{t} \nabla m(x,s)ds \big),
\end{split}
\end{align}
it follows that there exists $C >0$ such that
\begin{align}
\begin{split}
 |\nabla v(x,t)|
 \leq & |\nabla v_{0}(x)|+ \lambda v(x,t) |\int_{0}^{t} \nabla m(x,s) ds| \\
 \leq & C C_{0}^{2}
\end{split}
\end{align}
Since for all $(x,t) \in \O_T$
\begin{align}
\begin{split}
\Delta v(x,t) = & \Delta v_{0}(x) e^{-\lambda \int_{0}^{t} m(x,s)ds}
- 2 \lambda \nabla v_{0}(x) .  \big(\int_{0}^{t} \nabla m(x,s)ds \big) 
e^{-\lambda \int_{0}^{t} m(x,s)ds} \\
+ & \lambda^{2} v(x,t)
 \big| \int_{0}^{t} \nabla m(x,s)ds \big| ^{2}
- \lambda v(x,t) \big( \int_{0}^{t} \Delta m(x,s)ds \big)
\end{split}
\end{align}
it follows that
\begin{align}
\begin{split}
\label{smyi}
\forall (x,t) \in \O_T,\,\,\,\, |\Delta v(x,t)| \leq & C C_{0}^{3} +
 \lambda C_{0} |\int_{0}^{t} \Delta m(x,s)ds|
\end{split}
\end{align}
\noindent
with $C>0$ a suitable constant. \\
\noindent
For any fixed $x \in \O$, we integrate the equation $m_{t} - \alpha \Delta m + m = u$ on $[0,t]$
 and obtain that
$$\int_{0}^{t} \Delta m(x,s)ds =
\frac{1}{\alpha} [m(x,t) - m_{0}(x) + \int_{0}^{t} (m(x,s)-  u(x,s)) ds]$$
so that in view of \eqref{estm} there exists a constant $C>0$ such that
\begin{equation}
\label{smer}
\forall (x,t) \in \O_T,\,\,\,\, |\int_{0}^{t} \Delta m(x,s)ds|
\leq C C_{0}
\end{equation}
which completes the proof of $(b)$. Moreover in view of \eqref{smyi} and \eqref{smer},
 we conclude that there exists $C>0$ such that
$$\forall (x,t) \in \O_T,\,\,\,\, |\Delta v(x,t)| \leq C C_{0}^{3} $$
and obtain the property $(c)$, which completes the proof of Lemma \ref{xinjia}.

\subsection{Existence of a global solution to Problem $(P^{\varepsilon})$}
We prove the existence of a unique solution
$(u^{\varepsilon}, v^{\varepsilon}, m^{\varepsilon})$
to Problem $(P^{\varepsilon})$ on $\O_{T}$ for $\e>0$ small enough.
\begin{lem}
\label{new}
Assume that
$(u_{0}, v_{0}, m_{0})$ satisfy the hypotheses $1$-$2$-$3$-$4$. Then there exists
 $\varepsilon_{0} >0$ such that for all $0<\varepsilon<\varepsilon_{0}$,
Problem $(P^{\varepsilon})$ has a unique solution
 $(u^{\varepsilon}, v^{\varepsilon}, m^{\varepsilon})$
on $\O \times [0,T]$ for any $T>0$. This solution satisfies $0\leq \ue\leq C_0 $ in $\O_T$.
\end{lem}
The above lemma is similar to Lemma 4.2 in \cite{bhlm} and we just sketch the proof.
It relies on Schauder's fixed point theorem
and on the a priori estimates on Problem $(P^{\varepsilon})$ obtained in
Lemma \ref{xinjia}.

\noindent
First let $T>0$ be arbitrarily fixed and for all $u \in X_{T}$,
let $v= H(u)$ be defined as above. By the estimates of $v$ in Lemma \ref{xinjia},
there exists $ C  >0$
such that
\begin{equation}\label{vvv}
 0 \leq v \leq  C_0,\,\,\,
|\nabla v|  +
|\Delta v| \leq C C_0^3 \mbox{ in } \O_T.
\end{equation}
Next let $\tilde{u}$ be the unique solution of
\begin{equation}\label{uhat}
\left\{
\begin{array}{ll}
\tilde{u}_{t} = \Delta \tilde{u} - \nabla \cdot (\tilde{u} \nabla \chi (v)) + \dfrac{1}{\varepsilon^{2}} f(\tilde{u}) & \quad \quad \mbox{in} \ \O \times (0,T]  \\
\tilde{u}(x,0)=u_{0}(x) & \quad \quad x \in \Omega \\
\dfrac{\p \tilde{u}}{\p \nu}  =0  & \quad \quad \mbox{on} \ \p \O.
\end{array}
\right.
\end{equation}
The key point of the proof is to show that for $0<\e<\e_0$ small enough, we have
$$0 \leq {\hat u} \leq C_0 \mbox{ in } \O_T.$$
This follows from the fact that $C_0$ is a supersolution for equation (\ref{uhat})
for $\e>0$ small enough. Precisely,  using that $f(C_0)  <0$
since $C_0>1$ and  (\ref{vvv}),
we have that
\begin{eqnarray*}
 &&C_0 \Delta (\chi(v)) -
\frac{1}{\e^2} f(C_0)  \\
&=&C_0(\chi'(v)\Delta (v) + \chi''(v)|\nabla (v)|^2)-
\frac{1}{\e^2} f(C_0)\\
 &\geq&  - 2 C_0^4  -\frac{1}{\e^2} f(C_0)  \geq 0
\end{eqnarray*}
for $\e>0$ small enough. Moreover
$\tilde{u} \in C^{\alpha,\alpha/2}(\overline{\O_{T}})$ for some
$\alpha \in (0,1)$. Hence $u \rightarrow \tilde{u}$ maps $X_{T}$ into itself
and defines a compact operator.  A fixed point of this operator
obtained by Schauder's theorem
is then a solution
to Problem $(P^{\varepsilon})$. The uniqueness of solution follows from the a priori
 estimates on Problem $(P^{\varepsilon})$. For the details of the proof, we refer to
 \cite{bhlm} and \cite{dibe}.

\subsection{A comparison principle for Problem $(P^{\varepsilon})$}
We first recall the definition of a pair of sub- and super-solutions similar to the one
proposed in \cite{bhlm}.
\begin{df}
Let $(u_{\varepsilon}^{-}, u_{\varepsilon}^{+})$
be two smooth functions with
$ 0 \leq u_{\varepsilon}^{-} \leq u_{\varepsilon}^{+}$ in $\Omega_{T}$
and
$\dfrac{\p u_{\varepsilon}^{-}}{\p \nu} \leq \dfrac{\p u_{\varepsilon}^{+}}{\p \nu} $
 on $\p \O \times (0,T)$.
 By definition,
$(u_{\varepsilon}^{-}, u_{\varepsilon}^{+})$ is a pair of sub- and super-solutions in
 $\Omega_T$ if for any $v=H(u)$, with
$u_{\varepsilon}^{-} \leq u \leq u_{\varepsilon}^{+}$ in $\O_{T}$, we have
$$L_{v}[u_{\varepsilon}^{-}] \leq 0 \leq L_{v}[u_{\varepsilon}^{+}] \quad \mbox{  in  }
 \O_{T},$$
where the operator $L_{v}$ is defined by
$$L_{v}[\phi] = \phi_{t} - \Delta \phi + \nabla \cdot (\phi \nabla \chi (v) )
 - \frac {1}{\varepsilon^{2}} f(\phi).$$
\end{df}
\noindent
Note that in Lemma 2.2, $(0,C_{0})$ is a pair of sub- and super-solutions
 of Problem $(P^{\varepsilon})$.
It is then proved in \cite{bhlm} that the following comparison principle holds.
\begin{pro}
Let a pair of sub- and super-solutions $(u_{\varepsilon}^{-}, u_{\varepsilon}^{+})$
 in $\O_{T}$ be given. Assume that
$$\forall x \in \O,\,\,\,
u_{\varepsilon}^{-}(x,0) \leq u_{0}(x) \leq u_{\varepsilon}^{+}(x,0),$$
with $(u_{0}, v_{0}, m_{0})$ satisfying the hypotheses $1$-$2$.
Then there exists a unique solution $(u^{\varepsilon}, \ve, \me)$ of 
Problem $(P^{\varepsilon})$
with
$$\forall (x,t) \in \O_{T},\,\,\,
u_{\varepsilon}^{-}(x,t) \leq u^{\varepsilon}(x,t) \leq u_{\varepsilon}^{+}(x,t).$$
\end{pro}

\section{Well-posedness of Problem $(P^0)$}
\setcounter{equation}{0}
We establish here the existence and uniqueness of a smooth solution
 to the free boundary Problem $(P^{0})$ locally in time.
\begin{thm}
\label{well}
Let $\Gamma_{0} = \p \O_{0}$, where $\O_{0} \subset \subset \O$ is a $C^{2+\alpha}$
 domain with $\alpha \in (0,1)$.
Then there exists a time $T>0$ such that Problem $(P^{0})$
 has a unique solution $(v^{0},m^{0}, \Gamma)$ on $[0,T]$ with
$$\Gamma = (\Gamma_{t} \times \{ t \}) _{t \in [0,T]} \in C^{2+\alpha, (2+\alpha)/2}
  \mbox{   and   } v^{0}|_{\Gamma} \in C^{1+\alpha, (1+\alpha)/2}$$
\end{thm}
This theorem is similar to Theorem 2.1 in \cite{bhlm} and is using
a contraction fixed-point argument in suitable H\"older spaces
(see Section 2 in \cite{bhlm}).  We show here how it can actually be obtained using the
result established in
Theorem 2.1 in \cite{bhlm} and
some additional properties that we state
and prove below. \\
First we introduce some notations as in \cite{bhlm}.
We assume that $\Gamma_0$ is parametrized by some smooth $(N-1)$-dimensional
compact manifold ${\cal M}$ without boundaries which divides $\R^N$ into two pieces.
We denote by $ \vec{N}(s)$ the outward normal vector to ${\cal M}$
 at $s\in {\cal M}$
and define
$$\begin{array}{rll}
X:{\cal M}\times(-L,+L)& \rightarrow & \R^N \\
(s,s_N)&\mapsto& X(s,s_N)
\end{array}
$$
where
$$ X(s,s_N)= s + s_N \vec{N}(s). $$
If $L>0$ is chosen small enough, $X$ is a $C^{\infty}$-diffeomorphism
from ${\cal M}\times(-L,+L)$ onto a tubular neighborhood of ${\cal M}$ that
 we denote by ${\cal M}^L$.
We assume that $\Gamma_0 \subset {\cal M}^{L \over 2}$
 and is given by
$$\Gamma_0=\{ X(s,s_N),\, s_N=\Lambda_0(s),\, s\in {\cal M} \}$$
and that $\Omega_0$ is the connected component of $\Omega \setminus \Gamma_0$
 which contains
$$\{ x =X(s,s_N),\,  s_N < \Lambda_0(s), \,
s \in {\cal M} \}.$$
According to the regularity hypothesis on $\Gamma_0$ in 
Theorem \ref{well}, $\Lambda_0$ is a $C^{2 + \alpha}$ function with
$$||\Lambda_0||_{C^0({\cal M})} < \frac{L}{2}.$$
Let $T>0$ be a fixed constant that will be chosen later. 
We parametrize the interface $\Gamma= (\Gamma_t)_{t \in [0,T]}$ as follows
\begin{equation}
\label{gammat}
\Gamma_t=\{ X(s,s_N),\, s_N=\Lambda(s,t),\, s\in {\cal M} \},
\end{equation}
where $\Lambda: {\cal M} \times [0,T] \rightarrow (-L,+L)$ is a function.
 By definition, we will say that $\Gamma$ is $C^{m+\alpha,{m+\alpha \over 2}}$
 if the function $\Lambda$ satisfies
$$ \Lambda \in C^{m+\alpha,{m+\alpha \over 2}}({\cal M}\times [0,T]) $$
For any function $v(x,t)$ defined in $\overline{\O_T}$,
we consider the restriction of $v$ and of $\nabla v$ on the interface $\Gamma$ and
 we associate to $v$ the functions $w(s,t)$ and $\vec h(s,t)$ defined on
 ${\cal M}\times [0,T]$ by
\begin{eqnarray}
w(s,t)= v(X(s, \Lambda(s,t)), t), \label{w}\\
{\vec h}(s,t)= {\nabla v}(X(s, \Lambda(s,t)), t) . \label{h}
\end{eqnarray}
Next we split Problem $(P^{0})$ into two subproblems $(p_{a})$ and $(p_{b})$,
 where Problem $(p_{a})$ is given by
\begin{equation}
 (p_{a})
 \left\{
\begin{array}{ll}
V_{n} = -(N-1) \kappa + \chi^{'}(w){\vec h} \cdot {\vec n}
 \mbox{  on  } \ \Gamma_{t}= \p \O_{t},\,\, t \in (0, T]  \\
\Gamma_{t} |_{t=0} = \Gamma_{0}
\end{array}
\right.
\end{equation}
and Problem $(p_{b})$ is given by
\begin{equation}
 (p_{b})
 \left\{
\begin{array}{ll}
v_{t}^{0} =-\lambda m^{0}v^{0} &  \mbox{in}  \ \O \times (0,T] \\
m_{t}^{0}-\alpha \Delta m^{0} +  m^{0}= u^{0} &  \mbox{in}  \ \O \times (0,T] \\
\dfrac{\p m^{0}}{\p \nu} =0  &  \mbox{on} \ \p \O \times (0,T]\\
u^{0}(x,t)=  \chi_{\O_{t}}(x) = \displaystyle{
\left\{
\begin{array}{ll}
1  \mbox{   in } \ \O_{t}, t \in [0,T]   \\
0  \mbox{   in } \ \O \setminus \overline \O_{t}, t \in [0,T]
\end{array}
\right.}
\end{array}
\right.
\end{equation}
Note that the difference between the free boundary problem in \cite{bhlm} and here
concerns Problem $(p_{b})$.
Let us consider
\begin{equation}
 \label{M}
\forall (x,t) \in \Omega_T,\,\,\, M(x,t)= \int^{t}_{0} m^{0}(x,s) ds
\end{equation}
The restrictions
of $M$ and $\nabla M$ on $\Gamma$ are denoted
$a(s,t)$ and ${\vec b}(s,t)$ and defined on ${\cal M} \times [0,T]$ by
\begin{eqnarray}
a(s,t)= M(X(s, \Lambda(s,t)), t), \label{w}\\
{\vec b}(s,t)= {\nabla M}(X(s, \Lambda(s,t)), t) . \label{hb}
\end{eqnarray}
Note that using \eqref{guanyuv} and \eqref{erdba} we have that
$$w(s,t)=v_{0}(X(s, \Lambda(s,t))) e^{-\lambda a(s,t)}$$
and $${\vec h}(s,t)= \nabla v_{0}(X(s, \Lambda(s,t))) e^{-\lambda a(s,t)} -\lambda w(s,t) {\vec b}(s,t), $$
so that $w$ has the same regularity as $a$ and ${\vec h}$ has the same
 regularity as ${\vec b}$.
We deduce from Problem $(p_{b})$ that $M$ satisfies
\begin{equation}
 \label{youyige}
 \left\{
\begin{array}{ll}
-\alpha \Delta M+M=g(x,t) & \quad \quad \mbox{in}  \ \O \times (0,T] \\
\dfrac{\p M}{\p \nu} =0  & \quad \quad \mbox{on} \ \p \O \times (0,T],
\end{array}
\right.
\end{equation}
where $$g(x,t)=\int^{t}_{0} u^{0}(x,s)ds + m_{0}(x)- m^{0}(x,t). $$
 The same problem \ref{youyige} has been considered in \cite{bhlm}
but with a right-hand-side $g= u^0$.
Here the  function $g(x,t)$ is continuous in time, its regularity being the one of a time-integral of
$u^0$. Thus we can use Theorem 2.2 in \cite{bhlm}
and obtain (at least) the same regularity for $(a, {\vec b})$ 
in the case considered here.
\begin{lem}\label{t1}
Let $\Gamma=(\Gamma_t\times \{t\})_{t \in [0,T]}$ be given by \eqref{gammat} with
$$\Lambda \in C^{m + \alpha,{m + \alpha \over 2}}({\cal M}\times [0,T]) $$
for some $m \in \N$, $m \geq 2$ and $\alpha \in (0,1)$.
Let $M$ satisfy \eqref{youyige} and let $a$ and ${\vec b}$
 be associated to $M$ by (\ref{w}) and (\ref{hb}) respectively.
Then
$$a \in C^{m + \alpha,{m + \alpha \over 2}}({\cal M}\times [0,T]) $$
and
$${\vec b} \in [C^{m + \alpha',{m + \alpha'\over 2}}({\cal M}\times [0,T])]^n 
\mbox{ for all } 0<\alpha'< \alpha. $$
\end{lem}
 By the argument in \cite{bhlm} we know then that Problem $(p_{a})$ defines a mapping
 $(w,{\vec h}) \rightarrow \Lambda$ and Problem $(p_{b})$
 defines a mapping $\Lambda \rightarrow (w,{\vec h})$ with the proper 
regularity in H\"older
spaces. Therefore the composition of these two
mappings defines a contraction in some closed ball for $T>0$ small enough.
The unique fixed point of this contraction is the solution
to Problem $(P^{0})$ on $[0, T]$. This completes the proof of
Theorem \ref{well}.

\section{Generation of interface}
\setcounter{equation}{0}
In this section we establish the rapid formation of transition layers in a neighborhood
 of $\Gamma_{0}$ within a very short time interval of
order $\varepsilon^{2} |\ln \varepsilon|$. The width of the transition layer 
around $\Gamma_{0}$ is
 of order $\varepsilon$.
After a short time the solution $u^{\varepsilon}$ becomes
 close to 1 or 0 
except in a small
neighborhood of $\Gamma_{0}$. It reads precisely as follows.
\begin{thm}
\label{zao}
Let $u_{0}$ satisfy the assumptions $1$-$2$-$3$-$4$.
Let $0 <\eta < 1/4$ and define $\mu = f^{'}(1/2) = 1/4$.
Then there exist $\varepsilon_{0} >0$ and $M_{0}>0$ such that,
 for all $\varepsilon \in (0, \varepsilon_{0}]$ and $t^{\ast}=\mu^{-1} 
\varepsilon^{2} |\ln \varepsilon|$, \\
(a) for all $x \in \O$, we have
$$-\eta \leq u^{\varepsilon} (x, t^{\ast}) \leq 1+\eta ;$$
(b) for all $x \in \O$ such that $|u_{0}(x) - \frac{1}{2}| \geq M_{0} \varepsilon$, 
we have
$$\mbox{ if } \quad u_{0}(x) \geq \frac{1}{2} + M_{0} \varepsilon ,
  \mbox{ then }  u^{\varepsilon}(x, t^{\ast}) \geq 1-\eta ,$$
$$\mbox{ if } u_{0}(x) \leq \frac{1}{2} - M_{0} \varepsilon , \quad \mbox{then}  \quad u^{\varepsilon}(x, t^{\ast}) \leq \eta .$$
\end{thm}
\noindent
The above theorem relies on the construction of a suitable pair of 
sub- and super-solutions involving the solution of the bistable ODE.
We refer to the proof of Theorem 3.1 in \cite{a} 
in the simple case $\delta =0$.

\section{Convergence}
\setcounter{equation}{0}
We split the present section into 2 parts.
In a first step we establish the convergence
 of the $\ue$ to $u^0$ and prove  Corollary \ref{coru}). 
In a second step we prove Theorem \ref{thm1} as well as Theorem  \ref{gmi}, Theorem \ref{thm2} and 
Corollary \ref{haus}.

\noindent
In what follows, we construct a pair of sub- and super-solution
 $u_{\varepsilon}^{\pm}$ for Problem $(P^{\varepsilon})$ in order to control
 the function $u^{\varepsilon}$ on $[t^{\ast}, T]$. 
By the comparison principle it then follows that, 
if $u_{\varepsilon}^{-}(x,0) \leq u^{\varepsilon}(x,t^{\ast})
 \leq u_{\varepsilon}^{+}(x,0)$, 
then 
$u_{\varepsilon}^{-}(x,t) \leq u^{\varepsilon}(x,t+t^{\ast}) \leq u_{\varepsilon}^{+}(x,t)$ for all $(x,t) \in \O_{T}$. As a result, if both $u_{\varepsilon}^{+}$ and $u_{\varepsilon}^{-}$ converge to $u^{0}$, the solution $u^{\varepsilon}$ also converge to $u^{0}$ for all $(x,t) \in \O_{T} \setminus \Gamma$.

\subsection{Construction of sub- and super-solutions}
Before the construction, we present the definition of the modified 
signed distance function which is essential for our construction of 
sub- and super-solutions. Let us first define the signed distance function.
\begin{df}
Let $\Gamma = \bigcup _{0 \leq t \leq T} (\Gamma_{t} \times {t})$ be the solution of the limit geometric motion Problem $(P^{0})$. The signed distance function $\tilde{d}(x,t)$ is defined by
\begin{eqnarray}
&& \tilde{d}(x,t)= \displaystyle{
\left\{
\begin{array}{ll}
 dist(x, \Gamma_{t}) & \mbox{for} ~x \in \O \setminus \Omega_{t} \\
-dist(x, \Gamma_{t}) & \mbox{for}~ x \in \Omega_{t},
\end{array}
\right.}
\end{eqnarray}
\noindent
where $dist(x,\Gamma_{t})$ is the distance from $x$ to the hyperface $\Gamma_{t}$ in $\O$. \\
\noindent
Note that $\tilde{d}(x,t) =0$ on $\Gamma$ and that $|\nabla\tilde{d}(x,t)|=1$ in a neighborhood of $\Gamma$.
\end{df}
\noindent
In fact, rather than working with the above signed distance function $\tilde{d}(x,t)$, we need a modified signed distance function $d$ defined as follows.
\begin{df}
Let $d_{0} >0$ small enough such that $\tilde{d}(x,t)$ is smooth in
$$\{(x,t) \in \overline \O \times [0,T] , |\tilde{d}(x,t)| < 3 d_{0} \}$$
and such that for all $ t \in [0,T]$,
$$dist(\Gamma_{t}, \p \O) > 4 d_{0}.$$
\noindent
We define the modified signed distance function $d(x,t)$ by
$$d(x,t)= \zeta (\tilde{d}(x,t)),$$
where $\zeta (s)$ is a smooth increasing function on $\mathbb{R}^{N}$ defined by
\begin{eqnarray}
&& \zeta(s)= \displaystyle{
\left\{
\begin{array}{ll}
 s & \mbox{if} ~~|s| \leq 2d_{0} \\
-3d_{0} & \mbox{if} ~~ s \leq -3d_{0} \\
3d_{0} & \mbox{if} ~~ s \geq 3d_{0}.
\end{array}
\right.}
\end{eqnarray}
\end{df}
\noindent
Note that $|\nabla d|=1$ in the region 
$\{|d(x,t)| < 2d_{0} , (x,t) \in \overline \O \times [0,T] \}$. 
It follows that at $x \in \Gamma_t$,
the exterior normal vector  is $n(x,t) = \nabla d(x,t)$,
 the normal velocity is $V_n(x,t)= -d_t (x,t)$ and the mean curvature
is $K = \frac{1}{N-1}\Delta d(x,t)$. Therefore the motion law on $\Gamma^t$ given
 by Problem $(P^0)$ reads
\begin{equation}\label{eqd}
 d_t - \Delta d + \nabla d . \nabla \chi(v^0)=0  
\mbox{ on }
\Gamma_t=\{ x \in \Omega \Bigm|d(x,t)=0 \}.
\end{equation}
By Theorem \ref{gammat},  the interface $\Gamma_{t}$ is of class
 $C^{2+\alpha, \frac{2+\alpha}{2}}$ and $v^{0}$ is of class 
$C^{1+\alpha^{'}, \frac{1+\alpha^{'}}{2}}$ for any $\alpha, \alpha^{'} \in (0,1)$, 
all the functions $d_{t}$, $\Delta d$, $\nabla d$ are Lipschitz continuous
 near $\Gamma_{t}$ and $\nabla \chi(v^{0})$ is continuous near $\Gamma_{t}$. Therefore
from the mean value theorem applied separately on both sides of $\Gamma_{t}$, it follows
that there exists $N_0>0$ such that
\begin{equation}\label{eqdmv}
\forall (x,t) \in \Omega_T,\,\,\,
|d_t - \Delta d + \nabla d . \nabla \chi(v^0)| \leq N_0 |d(x,t)|.
\end{equation}
Note also that by construction, $\nabla d(x,t) =0$ in a neighborhood of $\p \O$. 

\noindent
As in \cite{a}, the sub- and super-solutions $u^{\pm} _{\varepsilon}$ are defined by
\begin{align}
\label{defyi}
\begin{split}
u_{\varepsilon}^{\pm}= U_{0}(\frac{d(x,t) \mp \varepsilon p(t)}{\varepsilon}) \pm q(t),
\end{split}
\end{align}
\noindent
where $U_{0}(z)$ is the unique solution of the stationary problem
\begin{equation}
\label{equyi}
\left\{
\begin{array}{ll}
U_{0}^{''} + f(U_{0}) =0  & \\
U_{0}(- \infty) =1 , U_{0}(0) = \frac{1}{2} , U_{0}( + \infty) = 0
\end{array}
\right.
\end{equation}	
\noindent
and
$$p(t) = -e^{-\beta t / \varepsilon^{2}} + e^{Lt} +K $$
$$q(t) = \sigma (\beta e^{-\beta t / \varepsilon^{2}} + \varepsilon^{2} L e^{Lt}) $$
with $L>0$ and $K>1$ to be chosen later. \\
\noindent
First note that $q = \varepsilon^{2} \sigma p_{t}$, then remark that for Problem \eqref{equyi} the unique solution $U_{0}$ has the following properties.
\begin{lem}
\label{lemyi}
There exist the positive constants $C$ and $\lambda$ such that the following estimates hold:
$$0 < U_{0}(z) \leq C e^{-\lambda |z|}  ~~\mbox{  for  } z \geq 0,$$
$$0 <1- U_{0}(z) \leq C e^{-\lambda |z|}  ~~\mbox{  for  } z \leq 0.$$
In addition, $U_{0} $ is strictly decreasing and $|U_{0}^{'}(z)| + |U_{0}^{''}(z)| \leq C e^{-\lambda |z|}$ for all $z \in \mathbb{R}$.
\end{lem}
The proof of Lemma 5.4 is given in \cite{bhlm}. We also note that
$$ u_{\varepsilon}^{-}(x,t) \leq  U_{0}(\frac{d(x,t)}{\varepsilon}) \leq u_{\varepsilon}^{+}(x,t)  $$
and that $p(t)$ is bounded for all $0 < \varepsilon < \varepsilon_{0}$
 and $t \in [0,T]$, $\lim_{\varepsilon \rightarrow 0} q(t) =0$ for all $t>0$. 
Therefore it follows from the definition of $u_\varepsilon ^\pm (x,t)$ that
for all $t \in (0,T]$,
\begin{eqnarray}
\label{uyl}
&& \lim_{\varepsilon \rightarrow 0} u_{\varepsilon}^{\pm}(x,t)=
 \chi_{\O_{t}}(x) =\displaystyle{
\left\{
\begin{array}{ll}
 1 & \mbox{for all} ~(x,t) \in \Omega_{t} \\
 0 & \mbox{for all} ~(x,t) \in \O \setminus \Omega_{t}
\end{array}
\right.}
\end{eqnarray}
\noindent
The key result of this section is the following lemma.
\begin{lem}
\label{lem1}
There exist $\beta>0,\sigma>0$ such that for all $K > 1$, 
we can find $\varepsilon_0>0$ and $L>0$ such that for any
 $\varepsilon \in (0,\varepsilon_0)$,  ($u_{\varepsilon}^{-}$, $u_{\varepsilon}^{+}$)
 is a pair of sub- and super-solutions for Problem $(P^{\varepsilon})$ 
in $\overline \O \times [0, T]$.
\end{lem}
\subsection{Proof of Lemma \ref{lem1}}
First note that for all $(x,t) \in \overline \O_{T}$,
$$u_{\varepsilon}^{-}(x,t) \leq U_{0}(\frac{d(x,t)}{\varepsilon})-q(t)
 \leq U_{0}(\frac{d(x,t)}{\varepsilon})+q(t) \leq u_{\varepsilon}^{+}(x,t).$$
Next since $\nabla d=0$ in a neighborhood of $\p \O$, we have that
 $\dfrac{\p u_{\varepsilon}^{\pm}}{\p \nu} = 0$ on $\p \O \times [0,T]$.
\noindent
Let $v$ be such that $v=H(u)$ with $u_{\varepsilon}^{-} \leq u 
\leq u_{\varepsilon}^{+}$ in $\O_{T}$, we show below that
$$L_{v}[u_{\varepsilon}^{-}] \leq 0 \leq L_{v}[u_{\varepsilon}^{+}],$$
where the operator $L_{v}$ is defined by
$$L_{v}[\phi] = \phi_{t} - \Delta \phi + \nabla(\phi \nabla \chi(v))
 - \frac {1}{\varepsilon^{2}} f(\phi).$$
Here we just consider the inequality $L_{v}[u_{\varepsilon}^{+}] \geq 0$,
 because the proof of the other inequality $L_{v}[u_{\varepsilon}^{-}] \leq 0$ 
is obtained by similar arguments.
A direct computation gives us the following terms
$$(u_{\varepsilon}^{+})_{t} = U_{0}^{'} (\frac{d_{t}}{\varepsilon} -p_{t}) +q_{t} ,$$
$$\nabla u_{\varepsilon} ^{+} =U_{0}^{'} \frac{\nabla d}{\varepsilon} ,$$
$$\Delta u_{\varepsilon} ^{+} = U_{0}^{''} \frac{|\nabla d|^{2}}{\varepsilon^{2}} 
+ U_{0}^{'} \frac{\Delta d}{\varepsilon} ,$$
where the value of the function $U_{
0}$ and its derivatives are taken at the point 
$\dfrac{d(x,t) - \varepsilon p(t)} {\varepsilon}$.
 Moreover the bistable function has the expansions
$$f(u_{\varepsilon}^{+}) = f(U_{0}) +q f^{'}(U_{0}) +
 \frac{1}{2} q^{2} f^{''}(\theta) , $$
where $\theta (x,t) $ is a function satisfying 
$U_{0} < \theta < u_{\varepsilon} ^{+}$. Hence, combining all the above,
 we obtain that
$$  L_{v}[u_{\varepsilon}^{+}]  =  
 (u_{\varepsilon}^{+})_{t} - \Delta u_{\varepsilon}^{+} 
+ \nabla u_{\varepsilon}^{+} \nabla \chi(v) 
+ u_{\varepsilon}^{+} \Delta \chi(v) -
\frac{1}{\varepsilon^{2}} f(u_{\varepsilon}^{+})
=   E_{1}+E_{2}+ E_{3}+E_{4}$$
where
$$E_{1} = -\frac{1}{\varepsilon^{2}} q[f^{'}(U_{0})
 + \frac{1}{2} q f^{''}(\theta)] - U_{0}^{'} p_{t} + q_{t},$$
$$E_{2} = \frac {U_{0}^{''}}{\varepsilon^{2}} (1- |\nabla d|^{2}),$$
$$E_{3} = \frac{U_{0}^{'}}{\varepsilon} (d_{t} - \Delta d 
+ \nabla d \cdot \nabla \chi(v_{0})),$$
$$E_{4} = \frac{U_{0}^{'}}{\varepsilon} 
\nabla d \cdot \nabla(\chi(v) - \chi(v^{0})) + u_{\varepsilon}^{+} \Delta \chi(v).$$
\noindent
We first need to present some useful inequalities before estimating 
the four terms above,
 this step is exactly the same as in \cite{a}. \\
Since $f^{'}(0) = f^{'}(1) = -\dfrac{1}{2}$ , 
we can find $0<b<1/2$ and $m>0$ such that
$$ \mbox{  if  } U_{0}(z) \in [0,b] \cup [1-b,1] 
\mbox{  then  } f^{'}(U_{0}(z)) \leq -m. $$
Furthermore, since the region $\{z \in \mathbb{R}, U_{0}(z) \in [b, 1-b] \}$ 
is compact and $U_{0}^{'} < 0$ on $\mathbb{R}$, there exists a constant $a_{1} >0$ 
such that
$$\mbox{  if  } U_{0} (z) \in [b,1-b] \mbox{  then  } U_{0}^{'}(z) \leq -a_{1}. $$
\noindent
Now we define $$F = \sup_{-1 \leq z \leq 2} (|f(z)| + |f^{'}(z)| + |f^{''}(z)|),$$
\begin{align}
\label{ok1}
\beta = \frac {m}{4},
\end{align}
and choose $\sigma$ which satisfies
\begin{align}
\label{ok2}
0 < \sigma < min (\sigma_{0}, \sigma_{1}, \sigma_{2}) ,
\end{align}
where $\sigma_{0} = \dfrac{a_{1}}{m+F}$, $\sigma_{1}= \dfrac{1}{\beta +1}$, $\sigma_{2} 
= \dfrac{4 \beta}{F(\beta+1)}$. Hence we obtain that
$$\forall z \in \mathbb{R}, -U_{0}^{'}(z) - \sigma 
f^{'}(U_{0}(z)) \geq 4 \sigma \beta.$$
\noindent
Now we have already chosen the appropriate $\beta$ and $\sigma$.
 Let $K>1$ be arbitrary, next we prove that $L_{v^{\varepsilon}}[u_{\varepsilon}^{+}]
 \geq 0$ provided that the constants $\varepsilon_{0}>0$ and $L>0$ 
are appropriately chosen. From now on, we suppose that the following inequality 
is satisfied
\begin{align}
\label{inyi}
\varepsilon_{0}^{2} L e^{LT} \leq 1.
\end{align}
\noindent
Then given any $\varepsilon \in (0, \varepsilon_{0})$, since $0<\sigma < \sigma_{1}$,
 we have $0 < q(t) < 1$ for all $t \geq 0$. Since $0 < U_{0} < 1$, 
it follows that for all $(x,t) \in \overline \O_{T}$
\begin{align}
\label{iner}
-1 < u_{\varepsilon}^{\pm}(x,t) < 2.
\end{align}
\noindent
We begin to estimate the four terms $E_{1}$, $E_{2}$, $E_{3}$ and $E_{4}$. 
The estimates of the terms $E_{1}$, $E_{2}$ and $E_{3}$ are similar
 to the estimates in \cite{a} and we obtain that
$$E_{1} \geq \frac{\sigma \beta^{2}}{\varepsilon^{2}} e^{-\beta t / \varepsilon^{2}} + 2 \sigma \beta L e^{Lt} =
\frac{C_{1}}{\varepsilon^{2}} e^{-\beta t / \varepsilon^{2}} + C_{1}^{'} L e^{Lt} ,$$
where $C_{1} = \sigma \beta^{2}$, $C_{1}^{'}= 2 \sigma \beta$ are positive constants.
$$|E_{2}| \leq 
\frac{16 C}{ (e \lambda d_{0})^{2}} (1+ ||\nabla d||_{\infty}^{2}) = C_{2} ,$$
where $C$ and $\lambda$ are the constants that we choose in Lemma 5.4, 
so that $C_{2}$ is also a positive constant. 

\noindent
We remark that in the estimate for $E_{2}$ in \cite{a}, the following assumption holds:
\begin{align}
\label{leq}
\begin{split}
e^{LT} + K \leq \frac{d_{0}}{2 \varepsilon_{0}} .
\end{split}
\end{align}
\noindent
For $E_{3}$, we use (\ref{eqdmv}) and obtain that
 $$|E_{3}| \leq C_{3} (e^{Lt} +K) + C_{3}^{'} ,$$
where $C_{3} =  N_{0} C$ and $C_{3}^{'} = \dfrac{N_{0} C}{\lambda}$
 with $C$ and $\lambda$ the constants given by Lemma 5.4 .

\noindent
Then we consider the term $E_{4}$. 
We should know the estimates
 of $\nabla (\chi(v)- \chi(v^{0}))$ and $\Delta \chi(v)$. 
In fact, for this term, we have the following lemma.
\begin{lem}
\label{lemer}
Let $u$ be any function satisfying 
 $$ u_{\varepsilon}^{-} \leq u \leq u_{\varepsilon}^{+} \mbox{ in } \O_{T}$$
and let $(v, m)$ be the corresponding solution of Problem \eqref{yeshi} with
 $v=H(u)$. Then there exists $C >0$ depending on $T$ and $\O$ such that
for all $(x,t) \in \O_{T}$, 
\begin{eqnarray}
&&|v(x,t)| + |\nabla v(x,t)| + |\Delta v(x,t)| \leq C \label{v2}\\
&&|\int_{0}^{t}(m-m^{0})(x,s)ds|+|\nabla d(x,t) \cdot \int_{0}^{t} \nabla (m-m^{0})(x,s)ds|
 \leq C \varepsilon p(t) \label{mdiff}\\
&&|(v-v^{0})(x,t)| + |\nabla d(x,t) \cdot \nabla(v-v^{0})(x,t)| \leq C
 \varepsilon p(t) \label{vdiff}
\end{eqnarray}
where $(v^{0},m^0)$ are given by the solution of Problem $(P^{0})$.
\end{lem}
\noindent
We prove this lemma below. Let us carry on with 
the proof of Lemma \ref{lem1}. We write
\begin{equation}
\label{wusisan}
\nabla d \cdot \nabla(\chi(v)-\chi(v^{0})) = \chi^{'}(v) \nabla d \cdot \nabla(v-v^{0}) + (\chi^{'}(v)-\chi^{'}(v^{0})) \nabla d \cdot \nabla v^{0}.
\end{equation}
Since $v^{0}$ is bounded in $C^{1+\alpha^{'}, \frac{1+\alpha^{'}}{2}}$
 for any $\alpha^{'} \in (0,1)$, 
there exists $C>0$, such that
$$ ||v^{0}||_{L^{\infty}(\O_{T})} +||\nabla v^{0}||_{L^{\infty}(\O_{T})} \leq C,  $$
which combined with \eqref{wusisan}, yields that
\begin{equation}
\label{wuwuling}
|\nabla d \cdot \nabla(\chi(v)-\chi(v^{0}))| \leq 
||\chi^{'}||_{\infty} |\nabla d \cdot \nabla(v-v^{0})| + C ||\nabla d||_{\infty} 
||\chi^{''}||_{\infty}|v-v^{0}|,
\end{equation}
where the $L^{\infty}$-norms of $\chi^{'}$ and $\chi^{''}$ 
are considered 
on the interval $(-C, C)$. Therefore, since $\chi$ is smooth 
and $||\nabla d||_{\infty} $ is bounded, 
it follows from \eqref{wuwuling} that for all $(x,t) \in \O_{T}$,
 there exists $C>0$ such that
\begin{equation}
\label{wuwuwu}
|\nabla d \cdot \nabla(\chi(v)-\chi(v^{0}))| \leq C \varepsilon p(t) .
\end{equation}
Moreover, using the smoothness of $\chi$ and the first inequality of Lemma 5.6,
 we obtain that there exists $C'>0$ such that
\begin{equation}
\label{wuliu}
|\Delta \chi(v)| \leq C'.
\end{equation}
Hence, by the above inequalities \eqref{wuwuwu}, \eqref{wuliu} 
and the fact that $| u_{\varepsilon}^{+}(x,t)| \leq 2$, we obtain that for all $(x,t) \in \O_{T}$,
$$|E_{4}| \leq \frac{C}{\varepsilon} C \varepsilon p(t)  + 2 C'. $$
\noindent
Finally substituting the expression for $p$ and $q$,
 we obtain that there exist the positive constants $C_{4}$, $C_{4}^{'}$ 
and $C_{4}^{''}$ such that
$$|E_{4}| \leq C_{4} + C_{4}^{'} e^{-\beta t / \varepsilon^{2}} + C_{4}^{''} e^{Lt}.$$
We collect the above four estimates of $E_{1}$, $E_{2}$, $E_{3}$ and $E_{4}$, which yield
\begin{align}
\begin{split}
L_{v} [u_{\varepsilon}^{+}]
\geq & \frac{C_{1}}{\varepsilon^{2}} e^{-\beta t / \varepsilon^{2}} 
+ C_{1}^{'} L e^{Lt} - C_{2} \\
         & - C_{3} (e^{Lt} +K) -    C_{3}^{'} - C_{4} 
- C_{4}^{'} e^{-\beta t / \varepsilon^{2}} 
- C_{4}^{''}  e^{Lt} \\
= & \frac{C_{1} - \varepsilon^{2} C_{4}^{'}}{\varepsilon^{2}} 
e^{-\beta t / \varepsilon^{2}} + (LC_{1}^{'} - C_{3} - 
C_{4}^{''}) e^{Lt} - C_{6} ,
\end{split}
\end{align}
where $C_{6} = C_{2} + C_{3}K + C_{3}^{'} + C_{4}$ is a positive constant.
\noindent
Now we set
$$L := \frac{1}{T} \ln \frac{d_{0}}{4 \varepsilon_{0}},$$
where $\varepsilon_{0}$ is small enough and satisfies the assumptions \eqref{inyi} and \eqref{leq}, so that $L$ is large enough. It also follows that $\dfrac{C_{1} - \varepsilon^{2} C_{4}^{'}}{\varepsilon^{2}} >0$ and
$$LC_{1}^{'} - C_{3} - C_{4}^{''} \geq \frac{1}{2} L C_{1}^{'},$$
therefore
$$L_{v} [u_{\varepsilon}^{+}] \geq \frac{1}{2} L C_{1}^{'} - C_{6} \geq 0.$$
The proof of Lemma 5.5 is now completed, with the constants $\beta$, $\sigma$ given in \eqref{ok1}, \eqref{ok2}.
\subsection{Proof of Lemma \ref{lemer}}
Lemma \ref{lemer} gives the key estimate and is the analogue of Lemma 4.9 in \cite{bhlm} 
and of Lemma 2.1 in \cite{a}. 
However the proof is markedly different since the coupling between $u$ and $v$ is given
 by a system with an ODE and a parabolic equation versus an elliptic equation in 
the two above references.

\noindent
First note that (\ref{v2}) is established exactly as in Lemma \ref{xinjia} (c).

\noindent
Concerning the second inequality (\ref{mdiff}), let us recall the following properties of $U_{0}$ 
given in [1].
\begin{lem}
\label{lemold}
For all given $a \in \mathbb{R}$ and $z \in \mathbb{R}$, we have the inequality:
$$|U_{0}(z+a) - \chi_{]-\infty,0]}(z)| \leq C e^{-\lambda |z+a|} + \chi_{]-a,a]}(z) $$
\end{lem}
Define $w(x,t) = m(x,t) - m^{0}(x,t)$, 
then $w$ satisfies
\begin{equation}
\label{side}
\left\{
\begin{array}{ll}
w_{t}-\alpha \Delta w + w = h  &\mbox{  in  } \O_{T}\\
\dfrac{\p w}{\p \nu}=0 & \mbox{  on  } \p\O \times (0,T) \\
w(x,0)=0, & x \in \O
\end{array}
\right.
\end{equation}
with $h = u - u^{0}$ satisfying
$$u_{\varepsilon}^{-} -u^{0} \leq \,h
 \leq \,  u_{\varepsilon}^{+} -u^{0} \mbox{ in  } \Omega_T.$$
 From the definition of $u_{\varepsilon}^{\pm}$ in \eqref{defyi} and from 
Lemma \ref{lemold}
 for $z= \dfrac{d(x,t)}{\varepsilon}$ and $a= \pm p(t)$, we deduce
 that for all $(x,t) \in \O_{T}$,
\begin{equation}
\label{ih}
|h(x,t)| \leq 
C(e^{- \lambda|d(x,t)/ \varepsilon + p(t)|} + e^{- \lambda|d(x,t)/ \varepsilon -p(t)|}) 
+ \chi_{\{ |d(x,t)| \leq \varepsilon p(t) \}} + q(t)
\end{equation}
Let us define for all $(x,t) \in \O_{T}$,
$$h_{1}(x,t)= q(t),$$
$$h_{2}(x,t)= C(e^{- \lambda|d(x,t)/ \varepsilon + p(t)|}
 + e^{- \lambda|d(x,t)/ \varepsilon -p(t)|}) 
 \chi_{\{ |d(x,t)| > d_{0} \}}$$
and
$$h_{3}(x,t)=C(e^{- \lambda|d(x,t)/ \varepsilon 
+ p(t)|} + e^{- \lambda|d(x,t)/ \varepsilon -p(t)|}) \chi_{\{ |d(x,t)| \leq d_{0} \}}
+ \chi_{\{ |d(x,t)| \leq \varepsilon p(t) \}} $$
and denote  by $(w_{i})_{i=1,2,3}$ the solutions 
of the three following auxiliary problems
$$
  (A_{i})
 \left\{
\begin{array}{ll}
(w_i)_t-\alpha \Delta w_i + w_i = h_i  &\mbox{  in  } \O_{T}\\
\dfrac{\p w_i}{\p \nu}=0 & \mbox{  on  } \p\O \times (0,T) \\
w_i(x,0)=0, & x \in \O
\end{array}
\right.
$$
Note that in view of the definition of $p(t)$ and the inequality \eqref{leq}, we have that for all $t \in [0,T]$
\begin{align}
\begin{split}
\label{ip}
0 < K-1 \leq p(t) \leq \frac{d_{0}}{2 \varepsilon_{0}}
\end{split}
\end{align}
so that the function $p$ is bounded away from $0$ 
for all $t \in [0,T]$. It follows in particular that choosing $\e>0$ small enough, 
$$ \varepsilon p(t)\leq d_{0}/ 2 \mbox{ for all }t \in [0,T]$$ so 
that $|h| \leq h_{1}+h_{2}+h_{3}$. Thus we deduce from the maximum principle 
that for all $x \in \O$ and $t \in [0,T]$, 
$$|w(x,t)| \leq w_{1}(x,t) + w_{2}(x,t)+w_{3}(x,t).$$
We now establish estimates for $w_i$, with $i=1,2,3$.

\noindent
\textbf{Problem $(A_{1})$} 

\noindent
Set $W_{1}(x,t)=\int_{0}^{t}w_1(x,s)ds$, then $W_{1}$ satisfies
\begin{equation}
\left\{
\begin{array}{ll}
(W_1)_t-\alpha \Delta W_{1} + W_{1} = H_{1}  &\mbox{  in  } \O_{T}\\
\dfrac{\p W_{1}}{\p \nu}=0 & \mbox{  on  } \p\O \times (0,T) \\
W_{1}(x,0)=0, & x \in \O
\end{array}
\right.
\end{equation}
with, since $q(t) = \varepsilon^{2} \sigma p'(t)$, 
$$H_{1}(x,t)=
\int_{0}^{t} q(s)ds= 
\varepsilon^{2} \sigma (p(t)-p(0))$$
so that by (\ref{ip}) we get that there exists 
$C>0$ such that for all $t \in [0,T]$,
$$\sup_{(y,s) \in \O \times [0,T]} |H_{1}(y,s)| \leq C \varepsilon p(t) .$$
Hence by standard parabolic estimates, there exists $C>0$ 
such that for all $(x,t) \in \O_{T}$,
 the solution $W_{1}$ of Problem $(A_{1})$ satisfies
\begin{equation}
\label{W1}
 |W_{1}(x,t)| + |\nabla W_{1}(x,t)| \leq C \varepsilon p(t).
\end{equation}
\textbf{Problem $(A_{2})$} 

\noindent
Note that by the standard parabolic estimates there exists a constant 
$C^{'}>0$ such that
By definition of $h_{2}$, using again \eqref{ip}, 
we obtain that there exists $C'>0$ such that for all $(s,t) \in [0,T]^{2}$
\begin{align}
\begin{split}
\label{liuliuyi}
h_{2}(y,s)  \leq & 2 C e^{-\lambda (d_{0}/ \varepsilon -p(s))} \\
\leq & 2C e^{-\lambda d_{0} / 2 \varepsilon}  \\
\leq & \frac{4C}{\lambda d_{0} e} \varepsilon  \\
\leq & \frac {4C}{ \lambda d_{0} e (K-1)} \varepsilon p(s)\\
\leq & C_{1} \varepsilon p(s) \leq  C'\varepsilon p(t) .
\end{split}
\end{align}
Thus by standard parabolic estimates, we obtain that for all $(x,t) \in \O_{T}$
$$|w_{2}(x,t)| + |\nabla w_{2}(x,t)| \leq C' \varepsilon p(t),$$
which implies that there exists $C>0$ such that for all $(x,t) \in \O_{T}$
\begin{equation}
 \label{W2}
|W_{2}(x,t)| + |\nabla W_{2}(x,t)| \leq C \varepsilon p(t),
\end{equation}
where we define $W_{2}(x,t)=\int_{0}^{t}w_2(x,s) ds$.

\noindent
\textbf{Problem $(A_{3})$} 

\noindent
Note that $h_{3}(y,s)$ is supported in $\{ |d(y,s)| \leq d_{0} \}$. 
Moreover by linearity 
we may suppose that the function $h_{3}$ satisfies
 one of the three following assumptions:
$$ (H_{1})  ~~~~  |h_{3}(y,s)| \leq \chi_{\{ |d(y,s)| \leq \varepsilon p(s) \}}$$
$$(H_{2}^{\pm})  ~~~~ |h_{3}(y,s)| \leq e^{-\lambda |d(y,s) / \varepsilon \pm p(s)|} $$
Then under respectively assumptions $(H_{1})$, $(H_{2}^{\pm})$, 
we define a function $\tilde{h}$ on $R \times [0,T]$, respectively by
\begin{eqnarray}
&& \tilde{h}(r,s)= \displaystyle{
\left\{
\begin{array}{ll}
 \chi_{\{ |r| \leq \varepsilon p(s) \}} \\
 e^{-\lambda |r/ \varepsilon \pm p(s)|}
\end{array}
\right.}
\end{eqnarray}
Note that $|h_{3}(y,s)| \leq \tilde{h} (d(y,s) ,s)$, 
and under either of the assumptions $(H_{1})$ or $(H_{2}^{\pm})$, 
there exists a constant $C >0$ such that for all $(s,t) \in [0,T]^{2}$
\begin{align}
\label{tildh}
\begin{split}
0 \leq \int_{-d_{0}}^{d_{0}} \tilde{h}(r,s)dr \leq C \varepsilon p(t).
\end{split}
\end{align}
Let $\varphi (x,t) = e^{t} w_{3}(x,t)$, then in view of Problem $(A_{3})$, the function
$\phi$ satisfies 
\begin{equation}
\label{cphi}
\left\{
\begin{array}{ll}
\varphi_{t}-\alpha \Delta \varphi= f  &\mbox{  in  } \O_{T}\\
\dfrac{\p \varphi}{\p \nu}=0 & \mbox{  on  } \p\O \times (0,T)
\end{array}
\right.
\end{equation}
where $f(x,t) = e^{t} h_{3}(x,t)$ and 
$\varphi(x,0) = w_{3}(x,0) = 0$ for all $x \in \O$. 
We establish now that 
there exist a constant $C>0$ such that 
\begin{equation}
 \forall (x,t) \in \Omega_T,\,\,\, 0\leq \varphi(x,t) \leq C \varepsilon p(t).
\label{fi}
\end{equation}
As in \cite{ahm}, the solution $\varphi(x,t)$ of Problem \eqref{cphi}
 can be expressed as
$$\varphi(x,t) = \int_{0}^{t} \int_{|d(y,s)| \leq d_{0}} G(x, y, t-s) f(y,s) dy ds, $$
with $G(x, y,t)$ 
being the Green function associated to the Neumann boundary value problem
 in $\O$ for the parabolic operator $\varphi_{t} -\alpha \Delta \varphi$.
Thus for all $(x,t) \in \Omega_T$,
\begin{equation}
 0 \leq \varphi(x,t)
\leq  
\int_{0}^{t} \int_{|d(y,s)| \leq d_{0}}
 G(x, y, t-s) e^s \tilde{h} (d(y,s) ,s) dy ds
\label{fir}
\end{equation}
Next we recall the following important property of $G$ which is established in
\cite{ahm}.

\noindent
{\bf Lemma 7.6, \cite{ahm}}: \textit{Let $\Gamma$ be a closed hypersurface
 in $\Omega$ and denote by $d(x)$ 
the signed distance function associated with $ \Gamma$. Then there exists constants
 $C, d_0>0$ such that for any function $\eta(r) \geq 0$ on $\R$, it holds that
$$\int_{|d| \leq d_{0}} G(x, y, t) \eta(d(y)) dy \leq 
\frac{C}{\sqrt{t}} \int_{-d_{0}}^{d_{0}} \eta(r)dr \mbox{ for } 0 < t \leq T $$}

\noindent
Moreover as pointed out in \cite{ahm}, the above inequality is uniform with respect to 
smooth variations of $\Gamma$ and for $t \in [0,T]$. 
Applying this inequality to our case, we deduce that there exists $C>0$ such that 
for all $(x,y) \in \Omega^2$
and for all $0\leq s < t \leq T$, 
\begin{equation}
 \label{intg}
\int_{|d(y,s)| \leq d_{0}} G(x, y, t-s) \tilde{h}(d(y,s),s) dy \leq 
\frac{C}{\sqrt{t-s}} \int_{-d_{0}}^{d_{0}} \tilde{h}(r,s) dr.
\end{equation}
In view of (\ref{fir}) and of (\ref{tildh}), it follows that for all $x \in \Omega$
and for all $t \in [0,T]$,  
\begin{eqnarray*}
 0\leq \varphi(x,t)
&\leq & 
C\int_{0}^{t} \int_{|d(y,s)| \leq d_{0}}
 G(x, y, t-s)  \tilde{h} (d(y,s) ,s) dy ds \\
&\leq  &
C'\int_{0}^{t}\frac{1}{\sqrt{t-s}}  \int_{-d_{0}}^{d_{0}} \tilde{h}(r,s) dr ds\\
&\leq  &
C'\int_{0}^{t}\frac{1}{\sqrt{t-s}}  \e p(t) ds \leq 2C' \e p(t) \sqrt{T}
\end{eqnarray*}
which yields inequality (\ref{fi}).

\noindent
Coming back to $w_3$, we deduce that for all $(x,t) \in \Omega_T$, 
\begin{equation}
\label{w3}
|w_{3}(x,t)| = |e^{-t} \varphi(x,t)| \leq C \varepsilon p(t).
\end{equation}
Define $W_{3}(x,t) = \int_{0}^{t} w_{3}(x,s)ds$, then it follows that
\begin{equation}
\label{W3}
|W_{3}(x,t)| \leq C \varepsilon p(t)
\end{equation}
 for some $C>0$ 
and for all $(x,t) \in \Omega_T$. 
We show now that there exist $C>0$ such that for all $(x,t) \in \O_{T}$,
\begin{equation}
\label{nablaW3}
|\nabla d(x,t).\nabla W_{3}(x,t)| \leq C \varepsilon p(t).
\end{equation}
Time integration of the equation in Problem $(A_3)$ on $[0,t]$ gives
$$w_{3}(x,t) - w_{3}(x,0) - \alpha \Delta W_{3}(x,t) +W_{3}(x,t) 
= \int_{0}^{t} h_{3}(x,s)ds.$$
 Since $w_3(x,0)=0$, 
we obtain the following elliptic problem for any $t \in [0,T]$,
\begin{equation}
\left\{
\begin{array}{ll}
-\alpha \Delta W_3(.,t) +W_3(.,t)=  \hat{H}_{3}(.,t) &\mbox{  in  } \O\\
\dfrac{\p W_{3}}{\p \nu}(.,t)=0 & \mbox{  on  } \p\O 
\end{array}
\right.
\end{equation}
where for all $(x,t) \in \O_{T}$,
$$\hat{H}_{3}(x,t) =  \int_{0}^{t} h_{3}(x,s)ds - w_{3}(x,t).$$
Let us define for any $t \in [0,T]$ the functions $a(.,t)$ as the solution of
\begin{equation}
\label{a}
\left\{
\begin{array}{ll}
-\alpha \Delta a(.,t) +a(.,t)=  h_3(.,t) &\mbox{  in  } \O\\
\dfrac{\p a}{\p \nu}(.,t)=0 & \mbox{  on  } \p\O 
\end{array}
\right.
\end{equation}
and define
$A(x,t) = \int_{0}^{t} a(x,s)ds$. Define similarly $B(.,t)$ as the solution of
\begin{equation}
\label{B}
\left\{
\begin{array}{ll}
-\alpha \Delta B(.,t) +B(.,t)=  -w_3(.,t) &\mbox{  in  } \O\\
\dfrac{\p B}{\p \nu}(.,t)=0 & \mbox{  on  } \p\O 
\end{array}
\right.
\end{equation}
so that by linearity 
$$\forall (x,t) \in  \Omega_T,\,\,\,W_3(x,t) =A(x,t) + B(x,t).$$
It follows from standard elliptic estimates in view of (\ref{w3}) that
$$|B(x,t)|+|\nabla B(x,t)| \leq C \varepsilon p(t).$$
Concerning $a$, note that the elliptic problem 
appearing here is the same as for the chemotaxis-growth system studied 
in \cite{bhlm} and in \cite{a}, with the right-hand-side satisfying
(\ref{tildh}).
Therefore the results stated in Lemma 4.2 in \cite{a} and in Lemma 4.10 in \cite{bhlm} 
 apply and prove that there exists a constant $C>0$ 
such that for all $(x,t) \in \O_{T}$,
$$|a(x,t)|+|\nabla d(x,t).\nabla a(x,t)| \leq C \varepsilon p(t)$$
and consequently
$$|A(x,t)|+|\nabla d(x,t). \nabla A(x,t)| \leq C \varepsilon p(t).$$
This completes
 the proof of (\ref{nablaW3}). 
In view of (\ref{W1}), (\ref{W2}), (\ref{W3}) and (\ref{nablaW3}),
inequality (\ref{mdiff}) is now established.

\noindent
In order to prove inequality (\ref{vdiff}), note that using (\ref{mdiff})
 we obtain that for all $(x,t) \in \O_{T}$, 
\begin{eqnarray}
|(v-v^{0})(x,t)| &= & |v_{0}(x) e^{-\lambda \int_{0}^{t} m(x,s)ds} 
- v_{0}(x) e^{-\lambda \int_{0}^{t} m^{0}(x,s)ds}| \\
&\leq & C|v_{0}(x)| |\int_{0}^{t} (m-m^{0})(x,s)ds| 
\leq  C' \varepsilon p(t),
\label{vv}
\end{eqnarray}
where $C'>0$ is a suitable constant. \\
Next we have similarly that for all $(x,t) \in \O_{T}$, 
\begin{eqnarray*}
&& |\nabla d(x,t) \cdot \nabla (v-v^{0})(x,t)| 
\leq C|e^{-\lambda \int_{0}^{t} m(x,s)ds}- 
e^{-\lambda \int_{0}^{t} m^{0}(x,s)ds}| \\
& &+ C |v(x,t) \nabla d(x,t)\cdot \int_{0}^{t} \nabla m(x,s)ds 
-v^{0}(x,t) \nabla d(x,t)\cdot \int_{0}^{t} \nabla m^{0}(x,s)ds| \\
&&\leq  C' \varepsilon p(t) + 
C |v(x,t)||\nabla d(x,t) \cdot \int_{0}^{t}  \nabla (m-m^{0})(x,s)ds| \\
& &+ C|v(x,t)-v^{0}(x,t)||\nabla d(x,t)
 \cdot \int_{0}^{t}
 \nabla m^{0}(x,s)ds |,
\end{eqnarray*}
where $C,C'>0$  are suitable constants. 
Using (\ref{vv}), (\ref{mdiff}) and upper bounds on $|v|$ and $|\nabla m^0|$,
we deduce that (\ref{vdiff}) is satisfied. This completes the proof of Lemma 5.6.
\subsection{Proof of Corollary \ref{coru} and Theorem \ref{thm1}}
The pointwise convergence of $u^{\varepsilon}$ to $u^{0}$ 
in $\bigcup_{0 < t \leq T} ((\O \setminus \Gamma_{t}) \times {t})$ 
when $\varepsilon \rightarrow 0$ follows from  Lemma \ref{lem1} and from (\ref{uyl}).
Next note that $w^{\e}=  \me - m^0$ is a solution of Problem \eqref{side}
 with the right-hand-side
$\he$ satisfying
$$ |\he(x,t)| \leq h_1(x,t) + h_2(x,t) + h_3(x,t)$$
with $h_i$, $i=1,2,3$ defined as in the proof of Lemma \ref{lemer}.
This shows that there exists  $C>0$ 
such that 
$$|\he||_{L^{1}(\O_T)} \leq C \varepsilon.$$
It follows then from standard parabolic estimates and Sobolev inequalities
 that for any $\alpha \in (0,1)$ there exists $p \in (1, +\infty]$ and  $C>0$ 
such that
\begin{align}
\begin{split}
||w^\e||_{C^{1+\alpha, 1+\alpha/2}(\overline \O_{T})} 
\leq & C ||u^{\varepsilon}-u^{0}||_{L^{p}(\O_{T})}\\
 \leq& C ||\he||_{L^{p}(\O_{T})} 
\leq   C \varepsilon^{\frac{1}{p}}.
\end{split}
\end{align}
Thus for any $\alpha \in (0,1)$,
$$\lim_{\varepsilon \rightarrow 0} 
||m^{\varepsilon}-m^{0}||_{C^{1+\alpha, (1+\alpha)/2}(\overline \O_{T})} =0.$$
The expression of $v^{\varepsilon}$ and $\nabla \ve$ in (\ref{guanyuv}) and (\ref{erdba})
then shows that 
$$\lim_{\varepsilon \rightarrow 0} 
||v^{\varepsilon}-v^{0}||_{C^{1+\alpha, (1+\alpha)/2}(\overline \O_{T})} =0$$
which completes the proof of Theorem \ref{thm1}.

\subsection{Proof of Theorem \ref{gmi}, Theorem \ref{thm2} and Corollary \ref{haus}}
The proofs are exactly the same as the proofs
 of Theorem 1.3, Theorem 1.5 and Corollary 1.6 in \cite{a} respectively,
 we omit the details here.


\end{document}